\documentclass[a4paper,10pt]{amsart}%
\usepackage{eurosym}
\usepackage{amsmath}
\usepackage{mathrsfs}
\usepackage{amsfonts}
\usepackage{graphicx}
\usepackage{amsfonts}
\usepackage{amssymb}%
\usepackage{esint}
\usepackage{color}
\usepackage[T1]{fontenc}
\usepackage[utf8]{inputenc}
\usepackage[backrefs]{amsrefs}
\usepackage{combelow}

\providecommand{\U}[1]{\protect\rule{.1in}{.1in}}
\newtheorem{theorem}{Theorem}[section]

\newtheorem{conjecture}[theorem]{Conjecture}
\newtheorem{openquestion}[theorem]{Open Question}

\newtheorem{proposition}[theorem]{Proposition}
\newtheorem{remark}[theorem]{Remark}

\newcommand{\eq}[1]{\begin{equation}{#1}\end{equation}}

\DeclareMathOperator{\M}{M}
\DeclareMathOperator{\I}{I}
\DeclareMathOperator{\diam}{diam}
\DeclareMathOperator{\supp}{supp}

\newcommand{\R}{\mathbb{R}}
\DeclareMathOperator{\BV}{BV}

\numberwithin{equation}{section}

\title{A trace inequality for solenoidal charges}

\author[B. Rai\cb{t}\u{a}]{Bogdan Rai\cb{t}\u{a}}
\author[D. Spector]{Daniel Spector}
\author[D. Stolyarov]{Dmitriy Stolyarov}

\address[B. Rai\cb{t}\u{a}]{Max-Planck-Instiutut f\"ur Mathematik in den Naturwissenschaften, Inselstrasse 22, 04103 Leipzig, Germany;\newline
Ennio De Giorgi Mathematical Research Center, Scuola Normale Superiore, Piazza dei Cavalieri 7, 56126 Pisa, Italy}
\email{bogdanraita@gmail.com}

\address[D. Spector]{Department of Mathematics, National Taiwan Normal University, No. 88, Section 4, Tingzhou Road, Wenshan District, Taipei City, Taiwan 116, R.O.C.;\newline
Okinawa Institute of Science and Technology Graduate University,
Nonlinear Analysis Unit, 1919--1 Tancha, Onna-son, Kunigami-gun,
Okinawa, Japan
}
\email{spectda@protonmail.com}

\address[D. Stolyarov]{St. Petersburg State University, Department of Mathematics and Computer Science,14th Line 29b, Vasilyevsky Island, St. Petersburg, Russia, 199178;\newline St. Petersburg Department of Steklov Mathematical Institute, Fontanka 27, St. Petersburg, Russia, 191023} 
\email{d.m.stolyarov@spbu.ru}

\begin{document}
\maketitle

\begin{abstract}
We prove that for $\alpha \in (d-1,d]$, one has the trace inequality
$$
\int_{\mathbb{R}^d} |\I_\alpha F| \;d\nu \leq C |F|(\mathbb{R}^d)\|\nu\|_{\mathcal M^{d-\alpha}(\R^d)}
$$
for all solenoidal vector measures $F$, i.e., $F\in M_b(\mathbb{R}^d,\mathbb{R}^d)$ and $\operatorname{div}F=0$. Here $\I_\alpha$ denotes the Riesz potential of order $\alpha$ and $\mathcal M^{d-\alpha}(\mathbb{R}^d)$ the Morrey space of $(d-\alpha)$-dimensional measures on $\mathbb{R}^d$.
\end{abstract}

\section{Introduction}
A result of V. Maz'ya \cite{mazya} (and later reproved by N. Meyers and W.P. Ziemer \cite{Meyers-Ziemer-1977}) asserts the existence of a constant $C_1>0$ such that one has the inequality
\begin{align}
\int_{\mathbb{R}^d} |u| \;d\nu \leq C_1\|\nu\|_{\mathcal{M}^{d-1}(\mathbb{R}^d)} |Du|(\mathbb{R}^d) \label{meyersziemer}
\end{align}
for every $u \in \BV(\mathbb{R}^d)$ and every non-negative Radon measure $\nu$ satisfying the ball growth condition $\|\nu\|_{\mathcal{M}^{d-1}(\mathbb{R}^d)}<+\infty$, where for $\beta \in (0,d]$,
\begin{align*}
\|\nu\|_{\mathcal{M}^{\beta}(\mathbb{R}^d)} :=   \sup_{x \in \mathbb{R}^d, r > 0} \frac{\nu(B(x,r))}{r^{\beta}}
\end{align*}
is the norm of $\nu$ in the Morrey space $\mathcal{M}^{\beta}(\mathbb{R}^d)$.  Here and in the sequel, with an abuse of notation, we use $|\cdot|$ to denote the total variation of vector-valued Radon measure (it also denotes the absolute value of a scalar and Euclidean norm of a vector).

The inequality \eqref{meyersziemer} is sometimes referred to as a trace inequality, as it gives an estimate for functions on lower dimensional subspaces, e.g. hyperplanes.  It is the analogue in the regime $p=1$ of the strong capacitary inequalities pioneered by V. Maz'ya \cite{Mazya2} and represents the state of the art concerning Sobolev inequalities under the assumption that $Du$ is a bounded measure.  Indeed, as discussed by A. Ponce and the second named author in \cite{Ponce-Spector}, and recounted in \cite[Section 6]{Spector-NA}, it implies the Sobolev inequality of E. Gagliardo \cite{Gagliardo} and L. Nirenberg \cite{Nirenberg}, its Lorentz improvement \cite{Alvino}, and even Hardy's inequality (the latter two are in fact equivalent in this case, as a result of the P\'olya-Szeg\"o inequality and \cite[Lemma 4.3 on p.~3424]{Frank-Seiringer}).  Yet perhaps the most significant consequence of the inequality \eqref{meyersziemer} are its implications concerning the well-definedness for $\mathcal{H}^{d-1}$ almost every $x \in \mathbb{R}^d$ not only of $u \in W^{1,1}(\mathbb{R}^d)$ but even $u \in \BV(\mathbb{R}^d)$, the deduction of which requires several additional ingredients found in the work of D. Adams \cite{Adams:1988}.  In particular, Adams shows firstly that the estimate for each such measure extends to an estimate on the space of functions Choquet integrable with respect to the Hausdorff content $\mathcal{H}^{d-1}_\infty$, $L^1(\mathcal{H}^{d-1}_\infty)$, and secondly that the Hardy--Littlewood maximal function is bounded on these spaces $L^1(\mathcal{H}^{\beta}_\infty)$, $\beta \in (0,d]$.  From this he obtains \cite[Theorem 1]{Adams:1988}:  There exists a constant $C_2=C_2(\alpha,d)>0$ such that
\begin{align}\label{adams}
\int_0^\infty  \mathcal{H}^{d-1}_\infty \left( \{ \M(u)>t\}\right)\;dt  \leq C_2|Du|(\mathbb{R}^d)
\end{align}
for all $u \in \BV(\mathbb{R}^d)$.  Here the integral on the left-hand-side is the Choquet integral with respect to the outer measure $\mathcal{H}^{d-1}_\infty$ and $\M$ is the Hardy-Littlewood maximal function, defined for $f \in L^1_{\mathrm{loc}}(\mathbb{R}^d)$ by
\begin{align*}
\M f(x):= \sup_{r>0}\fint_{B(x,r)} |f(y)|\;dy.
\end{align*}

In this paper we are interested in similar sharp trace inequalities for potentials acting on constrained subspaces of the space of vectorial measures, where surprisingly there are only two known results.  The first is not explicitly written, though it is almost an immediate consequence of \eqref{meyersziemer}: For $\alpha \in (1,d]$, there exists a constant $C_3=C_3(\alpha,d)>0$ such that
\begin{align}\label{curlfree}
\int_0^\infty  \mathcal{H}^{d-\alpha}_\infty \left( \{ \M (\I_\alpha F)>t\}\right)\;dt   \leq C_3 |F|(\mathbb{R}^d)
\end{align}
for all $F \in M_b(\mathbb{R}^d;\mathbb{R}^d)$ such that $\operatorname*{curl}F=0$ (or, equivalently, $F = D u$ for some~$u \in \dot{\BV}(\mathbb{R}^d)$); the symbols~$\I_\alpha$ and $M_b(\mathbb{R}^d;\mathbb{R}^d)$ denote the Riesz potential and  the space of~$\R^d$-valued measures (charges) of bounded total variation.  That the restriction $\alpha >1$ is necessary can be seen by the counterexample $F = D\chi_Q$, see e.g. \cite{Spector-PM}.  The second known result is \cite[Proposition 5]{Adams:1988}:  For $\alpha \in (0,d]$, there exists a constant $C_4=C_4(\alpha,d)>0$ such that
\begin{align}\label{hardy}
\int_0^\infty  \mathcal{H}^{d-\alpha}_\infty \left( \{ \M (\I_\alpha F)>t\}\right)\;dt  \leq C_4\|F\|_{\mathscr{H}^1(\mathbb{R}^d)}
\end{align}
for all $F$ in the real Hardy space $\mathscr{H}^1(\mathbb{R}^d)$.  In particular, the inequalities \eqref{curlfree} and \eqref{hardy} prompt one to wonder whether similar inequalities hold for other various constrained subspaces of measures, and if so, to determine the relationship between each subspace and the minimal $\alpha$ in such an inequality.  For curl free measures, it seems useful to express the validity of the inequality in terms of $\alpha>1=d-(d-1)$, since
$d-1$ is the largest number for which
\begin{align*}
\mathcal{H}^{d-1}(E) = 0 \implies |F|(E)=0
\end{align*}
for all $F \in M_b(\mathbb{R}^d;\mathbb{R}^d)$ with $\operatorname*{curl}F=0$.  This heuristic agrees with what one understands from the Hardy space inequality, where one has its validity for $\alpha>0=d-d$, $d$ being the largest number for which
\begin{align*}
\mathcal{H}^{d}(E) = 0 \implies |F|(E)=0
\end{align*}
for all $F \in \mathscr{H}^1(\mathbb{R}^d)$ (since any element of the Hardy space is absolutely continuous with respect to the Lebesgue measure).

These two examples motivate the following:  Given any closed translation and dilation invariant subspace $X \subset M_b(\mathbb{R}^d;\mathbb{R}^k)$, one defines a number associated with the dimension of the singular set of this subspace,
\begin{align*}
\kappa:= \inf_{\nu \in X} \operatorname*{dim_{\mathcal{H}}} \nu,
\end{align*}
where 
\begin{align*}
\operatorname*{dim_{\mathcal{H}}} \nu:= \sup_{s>0} \left \{ s : \mathcal{H}^{s}(E) = 0 \implies |\nu|(E)=0 \right \}.\end{align*}
Then one poses
\begin{openquestion}\label{qu:main}
	Let $\alpha \in (d-\kappa,d]$.  Can one show the existence of a constant $C=C(\alpha,d,X)>0$ such that
	\begin{align*}
	\int_0^\infty  \mathcal{H}^{d-\alpha}_\infty \left( \{\M(\I_\alpha F)>t\}\right)\;dt \leq C|F|(\mathbb{R}^d)
	\end{align*}
	for all $F \in X$?
\end{openquestion}

\begin{remark}
The answer to the open question above is positive in a related martingale model, see~\cite{ASW:2021}. It is also possible that one needs to impose some additional conditions on $X$ like closedness in some weaker topology, as in~\cite{Stolyarov:2020}. 
\end{remark}
 
The main result of this paper is to answer this question in the affirmative in the case of divergence free measures. Note that for such measures, it follows from Smirnov's theorem~\cite{Smirnov:1994} that $\kappa=1$.  In particular, we here establish
\begin{theorem}\label{TraceHausdorffContent}
Let $\alpha \in (d-1,d]$.  There exists a constant $C_5=C_5(\alpha,d)>0$ such that
\begin{align}\label{OurCapacitaryEstimate}
\int_0^\infty  \mathcal{H}^{d-\alpha}_\infty \left( \{ \M (\I_\alpha F)>t\}\right)\;dt \leq C_5|F|(\mathbb{R}^d)
\end{align}
for all $F \in M_b(\mathbb{R}^d;\mathbb{R}^d)$ such that $\operatorname*{div}F=0$.
\end{theorem}
\noindent

By the duality formula (see~\cite{Adams:1988} or Section 2.5 in~\cite{AdamsHedberg:1996})
\begin{align}\label{functional_equivalence}
\int_0^\infty \mathcal{H}^\beta_\infty\big(\{x\mid |g(x)| \geq t\}\big)\,dt \asymp \sup\Big(\Big\{\int\limits_{\mathbb{R}^d} g(x)\,d\mu(x)\,\Big|\; \|\mu\|_{\mathcal{M}^\beta} \leq 1\Big\}\Big),
\end{align}
and the boundedness of the maximal function on $L^1(\mathcal{H}^\beta_\infty)$, Theorem~\ref{TraceHausdorffContent} is equivalent to the trace inequality given in
\begin{theorem}\label{OurTrace}
Let $\alpha \in (d-1,d]$.  There exists a constant $C_6=C_6(\alpha,d)>0$ such that
\begin{equation}\label{OurTraceEstimate}
\int_{\mathbb{R}^d} |\I_\alpha F| \;d\nu \leq C_6\|\nu\|_{\mathcal{M}^{d-\alpha}(\mathbb{R}^d)} |F|(\mathbb{R}^d)
\end{equation}
for all vector measures~$F\in M_b(\R^d;\R^d)$ such that $\operatorname*{div}F=0$.
\end{theorem}

Theorems~\ref{TraceHausdorffContent} and \ref{OurTrace} are sharp in the sense that they fail for $\alpha \in (0,d-1]$.  Indeed, as in D. Adams' proof of \cite[Proposition 5 on p.~121]{Adams:1988}, the proof of any value in this range would imply the validity of the result for $\alpha=d-1$, which cannot hold as a result of the following
\begin{theorem}\label{rk:no_endpoint_est}
There exists $F \in M_b(\mathbb{R}^d;\mathbb{R}^d)$ with $\operatorname*{div}F=0$ in the sense of distributions and
\begin{align*}
\sup_{t>0} t\mathcal{H}^{1}_\infty\left(\{ |\I_{d-1} F|>t\}\right)  =+\infty.
\end{align*}
\end{theorem}

Inequalities analogous to ~\eqref{OurCapacitaryEstimate} and~\eqref{OurTraceEstimate} hold for broader classes of differential constraints.  In particular, if one has a measure $F$ and a first order cocancelling (see ~\cite{SpectorHernandezRaita2021,VanSchaftingen2013} for a definition, for example) differential constraint $L$ for which $L(D)F=0$, then one can write
\begin{align*}
F=T^\dagger T F
\end{align*}
where $T,T^\dagger$ are maps on finite dimensional spaces and $\operatorname*{div} TF=0$ row-wise.  Therefore, the estimate for divergence free fields extends to those which admit such annhilators.  We refer to \cite{SpectorHernandezRaita2021} for the details, which is based upon an idea from \cite{VanSchaftingen2013}.  Note, however, that for other choices of first order differential constraints the result may fail to be sharp with respect to the minimal admissible $\alpha$.  For example, our result implies the curl free case, inequality \eqref{curlfree}, for $\alpha>d-1$, while the result in fact holds for any $\alpha>1$.

Concerning the endpoint case $\alpha = d - 1$ (and more generally $\alpha=d-\kappa$), the situation remains unclear.  In the curl free case, while the inequality \eqref{curlfree} fails at the endpoint $\alpha=1$, the inequality \eqref{adams} is only a singular integral transformation away:  If we denote by $R^*$ the adjoint of the vector-valued Riesz transform, while the estimate fails for $\I_1 \nabla u$, it holds for $R^*\cdot \I_1 \nabla u=u$.  This suggests that with an appropriate singular integral transformation (with additional cancellation properties) it may be possible to obtain an estimate in the endpoint $\alpha=d-1$ (and by extension, in endpoints for other various subspaces). In particular, we let~$K_1,K_2,\ldots,K_d$ be a collection of (sufficiently smooth) functions on~$\R^d$ that are homogeneous of order~$-1$ and consider the operator
\eq{
\mu \mapsto \sum\limits_{j=1}^d K_j*\mu_j
}
acting on~$\R^d$-valued charges. Let us call this operator~$K$. 
\begin{conjecture}
The inequality
\eq{
\int_{\mathbb{R}^d} |K[F]| \;d\nu  \lesssim |F|(\mathbb{R}^d) \|\nu\|_{\mathcal{M}^{1}(\mathbb{R}^d)},\qquad \operatorname*{div} F = 0,
}
holds true if and only if 
\eq{\label{WeakCanc}
\sum\limits_{j=1}^d K_j(\xi)\xi_j = 0.
}
for any~$\xi \in \R^d$.
\end{conjecture}
The ``only if'' may be obtained by testing the case where~$F$ and~$\nu$ are concentrated on one and the same segment. The condition~\eqref{WeakCanc} may be restated in terms of the Fourier transform as~$\mathrm{div}[\hat{\mathcal{K}}] = 0$, where~$\mathcal{K} = (K_1,K_2,\ldots,K_d)$.

\section{Proofs}\label{proofs}
The notation $A\lesssim B$ means there exists~$C > 0$ such that~$A \leq CB$, where $C$ may depend on the dimension, $\alpha$, but not on the functions or measures being estimates.  For example, in formula~\eqref{Morrey1} it does not depend on the choice of~$\mu$. The following technical and elementary proposition is interesting in itself. 
\begin{proposition}\label{AdamsHedberg}
Let~$\mu$ be a \textup(signed or vector valued\textup) measure on~$\R^d$ with compact support, zero mean, and such that
\eq{\label{Morrey1}
\|\mu\|_{\mathcal{M}^1} \lesssim 1.
}
Then, for any~$\alpha \in (d-1,d)$ the inequality
\eq{
\int_{\mathbb{R}^d} |\I_\alpha \mu| \;d\nu \lesssim \diam (\supp \mu)\cdot\|\nu\|_{\mathcal{M}^{d-\alpha}}
}
holds true for any non-negative measure~$\nu$.
\end{proposition}
\begin{proof}
Without loss of generality, by translation we may assume $0 \in \supp \mu$.  Define
\begin{align*}
R:= \max_{x \in \supp \mu} |x|,
\end{align*}
so that $\mu$ is supported in the ball~$B(0,R)$ with $R\leq 2\diam (\supp \mu)$. We will estimate~$\I_\alpha\mu$ at the points~$x\in B(0,2R)$ and~$x\notin B(0,2R)$ in two different ways. 

Let us start with the former case, where a telescoping dyadic argument and the location of $x$ in relation to the support of $\mu$ yields the inequality
\begin{align*}
\big|\I_{\alpha}\mu(x)\big| \lesssim \int\limits_{\R^d} \frac{d|\mu|(y)}{|x-y|^{d-\alpha}} \lesssim \sum\limits_{2^k \leq 3R} 2^{-(d-\alpha-1)k}\|\mu\|_{\mathcal{M}^1} \lesssim R^{-d+\alpha +1}.
\end{align*}

Concerning the latter case, using that~$\int \,d\mu = 0$ we have
\begin{align*}
|\I_{\alpha}\mu(x)|  = c_{d,\alpha}\Big|\int\limits_{\R^d} \frac{d\mu(y)}{|x-y|^{d-\alpha}}\Big| = c_{d,\alpha}\Big|\int\limits_{\R^d} \Big(\frac{1}{|x-y|^{d-\alpha}} - \frac{1}{|x|^{d-\alpha}}\Big)d\mu(y)\Big|.
\end{align*}
In particular, an application of the mean value theorem gives the standard estimate
\begin{align*}
\Big|\frac{1}{|x-y|^{d-\alpha}} - \frac{1}{|x|^{d-\alpha}}\Big| \lesssim \frac{|y|}{|x|^{d-\alpha + 1}},\qquad y\in B(0,R),\ x\notin B(0,2R),
\end{align*}
so that in this regime
\begin{align*}
|\I_{\alpha}\mu(x)| \lesssim |x|^{-d+\alpha - 1}\int\limits_{\R^d} |y|\,d|\mu|(y)  \lesssim |x|^{-d+\alpha-1} R^2.
\end{align*}

Therefore, we have proved the estimate
\eq{
\Big|\I_{\alpha}\mu(x)\Big| \lesssim 
\begin{cases}
R^{-d+\alpha + 1},\qquad & x\in B(0,2R); \\ 
|x|^{-d+\alpha-1} R^2,\qquad & x\notin B(0,2R).
\end{cases}
}
We integrate this estimate with respect to~$\nu$:
\begin{align*}
\int\limits_{\R^d}|\I_\alpha \mu(x)|\,d\nu(x) &\lesssim \int\limits_{B_{2R}(0)} R^{-d + \alpha + 1}\,d\nu(x) + R^2\int\limits_{\R^d \setminus B_{2R}(0)} |x|^{-d+\alpha -1}\,d\nu(x) \\
 &\lesssim R\|\nu\|_{\mathcal{M}^{d-\alpha}} + \sum\limits_{2^k \geq R}2^{k(-d+\alpha-1)} 2^{k(d-\alpha)} R^2\|\nu\|_{\mathcal{M}^{d-\alpha}} \\
 &\lesssim R\|\nu\|_{\mathcal{M}^{d-\alpha}}.
\end{align*}
\end{proof}
We next give the
\begin{proof}[Proof of Theorem~\ref{OurTrace}]
The case $\alpha=d$ is an $L^\infty$ estimate which has been established in the literature (e.g. see Theorem~$5$ in~\cite{Raita}), therefore we focus on the case $\alpha \in (d-1,d)$.  Following the argument in \cite{SpectorHernandez2020} and \cite{SpectorHernandezRaita2021}, by Smirnov's decomposition (see~\cite{Smirnov:1994}) and the surgery lemma (Lemma~$4.1$ in~\cite{SpectorHernandez2020}), it suffices to prove the estimate for $F=\mu_\Gamma$, where
$\mu_\Gamma$ is a measure induced by integration along a piecewise-$C^1$ closed loop $\Gamma$:
\begin{align*}
\int_{\mathbb{R}^d} \Phi\cdot \mu_\Gamma := \int_0^{|\Gamma|} \Phi(\gamma(t))\cdot \dot{\gamma}(t)\;dt,\qquad \Phi \in C(\R^d,\R^d),
\end{align*}
and which satisfies
\begin{align*}
\|\mu_\Gamma\|_{\mathcal{M}^{1}} \lesssim 1.
\end{align*}  
Here we use the notation $\gamma:[0,|\Gamma|] \to \mathbb{R}^d$ to denote the parametrization of the closed loop $\Gamma$  by arclength.  

However, for any such curve we have that the diameter of the curve is proportional to its total variation, $\diam \supp \Gamma \lesssim \|\mu_\Gamma\|_{M_b(\mathbb{R}^d;\mathbb{R}^d)}$, and therefore Proposition~\ref{AdamsHedberg} implies the theorem.
\end{proof}

\begin{proof}[Proof of Theorem~\ref{rk:no_endpoint_est}]
Define the curve $\Gamma$ to be the boundary of the square $(0,1)^2$ embedded in $\R^2\times\R^{d-2}$ (which in the sequel we denote by $(x_1,x_2,x')$).  For such a choice of $\Gamma$, we let $\gamma:[0,4]\to \mathbb{R}^d$ be its parametrization by arclength, oriented counterclockwise.  Then the desired solenoidal measure is $F=\dot{\gamma}\mathcal{H}^1|_{\Gamma}$.  Indeed, one has that $F \in M_b(\mathbb{R}^d;\mathbb{R}^d)$ and is divergence free, the latter following from the fact that it is closed (one can see this by an application of the fundamental theorem of calculus and using that the endpoints are the same).  Meanwhile, $\I_{d-1}$ applied to $F$ is given by:
\begin{align*}
	\I_{d-1}(\dot{\gamma}\mathcal{H}^1|_{\Gamma})(x)=\int_{\Gamma} \dfrac{\dot{\gamma}(y)d\mathcal{H}^1(y)}{|x-y|} \quad\text{for }x\in\R^d.
\end{align*}

\begin{figure}
\includegraphics[width = 1\linewidth]{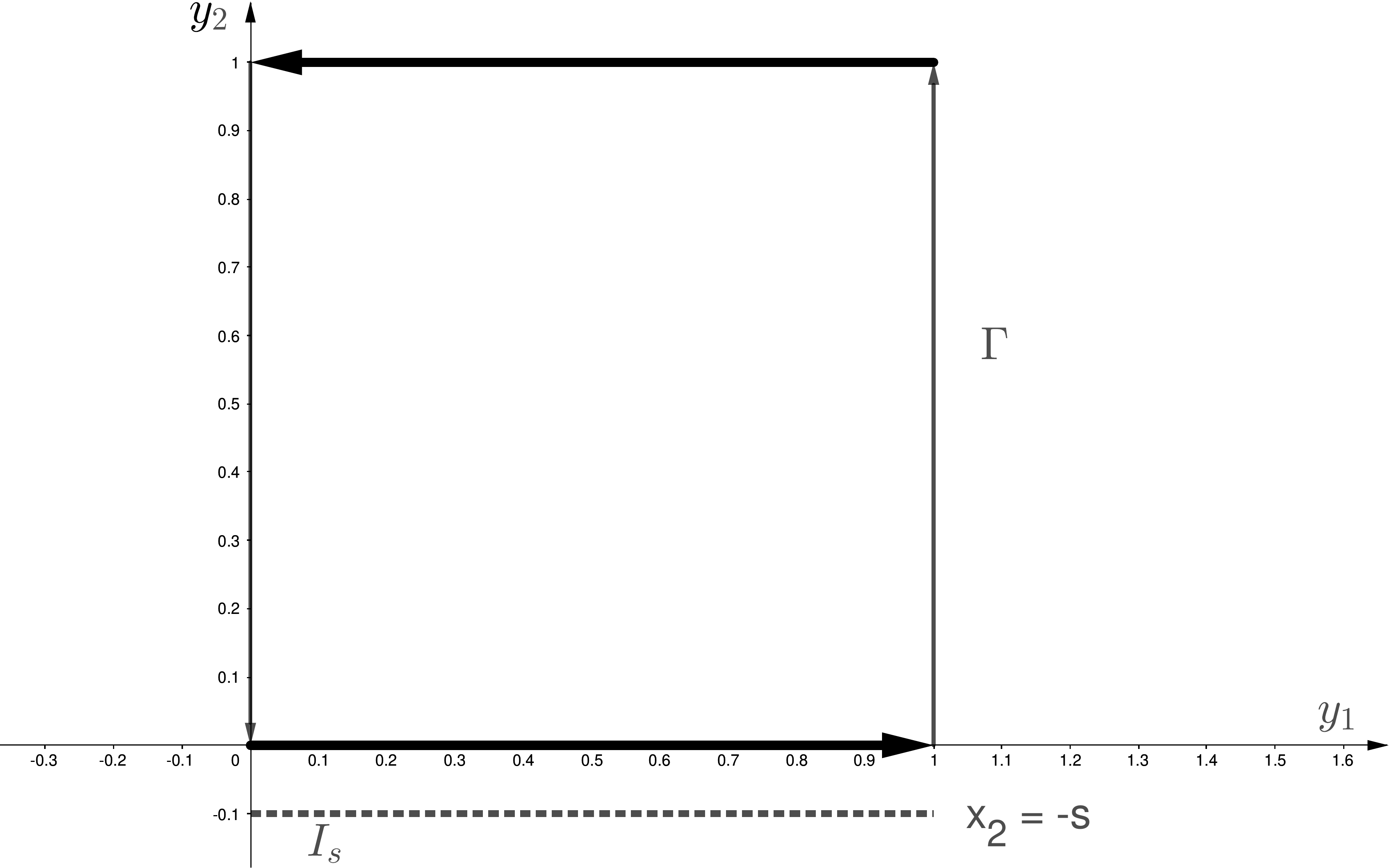}
\caption{Illustration to the proof of Theorem~\ref{rk:no_endpoint_est}.}
\end{figure}

For our purposes it will suffice to make estimates for the first component,
\begin{align}\label{LongFormula}
		[\I_{d-1}(\dot{\gamma}\mathcal{H}^1|_{\Gamma})(x)]_1&=\left(\int_{\Gamma\cap\{y_2=0\}}+\int_{\Gamma\cap\{y_2=1\}}\right) \dfrac{\dot{\gamma}_1(y)d\mathcal{H}^1(y)}{|x-y|}\\
&		= \int_0^1\dfrac{dy_1}{(|x_1-y_1|^2+|x_2|^2+|x'|^2)^{1/2}}\\
&-\int_0^1\dfrac{dy_1}{(|x_1-y_1|^2+|x_2-1|^2+|x'|^2)^{1/2}}.
\end{align}

In particular, to prove the failure of the weak-type estimate with respect to the content it suffices to show that the estimate blows up for the first component,
\begin{align*}
\sup_{t>0} t\mathcal{H}^1_\infty\left(\{ |\I_{d-1}(\dot\gamma_1\mathcal{H}^1|_{\Gamma})|>t\}\right)=\infty.
\end{align*}
By the functional equivalence given in \eqref{functional_equivalence}, this blow up, in turn, will be demonstrated if we can show that for any $t>0$ sufficiently large, one can find a measure
$\mu_s$, $s=s(t)$, with $\|\mu_s\|_{\mathcal{M}^1} \l \leq 1$ and \begin{align}\label{blowup}
\mu_s\left(\{|\I_{d-1}(\dot\gamma_1\mathcal{H}^1|_{\Gamma})|>t\}\right) \geq 1.
\end{align}

Returning to our choice of $\Gamma$, in the halfspace given by $x_2\leq 0$, the second term  in~\eqref{LongFormula} is bounded from above by $1$, so
\begin{align*}
		|\I_{d-1}(\dot{\gamma}_1\mathcal{H}^1|_{\Gamma})(x)|&\geq \int_0^1\dfrac{dy_1}{(|x_1-y_1|^2+|x_2|^2+|x'|^2)^{1/2}}-1.
\end{align*}
Thus, if $x_1 \in (0,1)$, $x_2=-s<0$ for some $s \in (0,1)$ and $x'=0$, since $\max\{1-x_1,x_1\}\geq 1/2$, we have
\begin{align*}
|\I_{d-1}(\dot{\gamma}_1\mathcal{H}^1|_{\Gamma})(x)|&\geq \int_{\frac{-x_1}{s}}^{\frac{1-x_1}{s}} \dfrac{dz}{(z^2+1^2)^{1/2}}-1\\
&\geq \int_0^{1/2s}\dfrac{dz}{(z^2+1)^{1/2}}-1\\
&=\ln \left(\sqrt{1+1/4s^2}+1/2s\right)-1\\
&\geq \ln 1/s-1.
\end{align*}

In particular, for every $s \in (0,1)$ we define the measures $\mu_s$ by $\mu_s=\mathcal{H}^1|_{I_s}$ where $I_s=\{(x_1,-s,0)\colon x_1\in(0,1)\}$.  Then $\|\mu_s\|_{\mathcal{M}^1} \l \leq 1$ and for $x \in \supp \mu_s$ one has
\begin{align*}
		|\I_{d-1}(\dot{\gamma}_1\mathcal{H}^1|_{\Gamma})(x)|&\geq \ln 1/s-1.
\end{align*}
Therefore
\begin{align*}
\mu_s\left(\{|\I_{d-1}(\dot\gamma_1\mathcal{H}^1|_{\Gamma})|>t\}\right) \geq \left|\left\{x_1 \in (0,1): \ln 1/s>t+1\right\}\right|.
\end{align*}
However, the condition on the right hand side is uniform over $x_1 \in (0,1)$, and for every $t>0$ sufficiently large, any choice of $s < \exp(-t-1)$ yields
\begin{align*}
\mu_s\left(\{|\I_{d-1}(\dot\gamma_1\mathcal{H}^1|_{\Gamma})|>t\}\right) \geq 1,
\end{align*}
which completes the proof of the claimed inequality \eqref{blowup} and therefore the Theorem.
\end{proof}

\section*{Acknowledgments}
D. Stolyarov is partially supported by RFBR grant  20-01-00209 and by ``Native towns'', a social investment program of PJSC ``Gazprom Neft''. Part of this work was undertaken while D. Spector was visiting the National Center for Theoretical Sciences in Taiwan.  He would like to thank the NCTS for its support and warm hospitality during the visit.

\end{document}